\documentclass{article}

\usepackage{amsmath,amssymb,theorem}

\newcommand{\al}{\alpha}
\newcommand{\bt}{\beta}

\newcommand{\ov}{\overline}
\newcommand{\C}{\mathbb C}
\newcommand{\F}{\mathbb F}
\newcommand{\Q}{\mathbb Q}
\newcommand{\R}{\mathbb R}
\newcommand{\Z}{\mathbb Z}

\newcommand{\Sym}{\operatorname{Sym}}

\newcommand{\mult}{\operatorname{mult}}
\newcommand{\sign}{\operatorname{sign}}

\newcommand{\eop}{\hspace*{\fill} $\Box$}

\theoremstyle{break}

\newtheorem{thm}{Theorem}[section]
\newtheorem{prp}[thm]{Proposition}

\begin{document}

\begin{center}
{\Large\bf Natural constructions of some generalized}\\[0.2cm] 
{\Large\bf Kac-Moody algebras as bosonic strings}\\[2cm]

Thomas Creutzig, Alexander Klauer, Nils R. Scheithauer\\[2cm]
\end{center}
There are 10 generalized Kac-Moody algebras whose denominator
identities are completely reflective automorphic products of singular
weight on lattices of squarefree level. Under the assumption that the
meromorphic vertex operator algebra of central charge 24 and spin-1
algebra $\hat{A}_{p-1,p}^r$ exists we show that four of them can be
constructed in a uniform way from bosonic strings moving on suitable
target spaces. 
\vspace*{1.5cm}

\section{Introduction}

Generalized Kac-Moody algebras are natural generali\-zations of the finite
dimensional simple Lie algebras. They are defined by generators and
relations and are allowed to have imaginary simple roots. Generalized
Kac-Moody algebras are in general infinite dimensional but their
theory is in many aspects similar to the finite dimensional theory. In
particular there is a character formula for highest weight modules and
a denominator identity. Borcherds has used twisted versions of the
denominator identity of the monster algebra to prove Conway and
Norton's moonshine conjecture \cite{B3}. 

Borcherds' singular theta correspondence \cite{B5} is a map from
modular forms for the Weil representation to automorphic forms on
orthogonal groups. Since these automorphic forms can be written as
infinite products they are called automorphic products. They have
found various applications in geometry, arithmetic and in the theory
of Lie algebras. In particular the denominator identities of generalized
Kac-Moody algebras are sometimes automorphic products (cf.\ \cite{S1},
\cite{S2}). Reflective automorphic products are automorphic products
whose di\-vi\-sors correspond to roots of the underlying lattice and
are zeros of order one. They can be classified under certain
conditions \cite{S3}. 

In \cite{S1} a family of 10 generalized Kac-Moody algebras is
constructed. Recall that the Mathieu group $M_{23}$ acts on the Leech
lattice $\Lambda$. Let $g$ be an element of squarefree order $N$ in
$M_{23}$. Then $g$ has characteristic polynomial \mbox{$\prod_{k|N}
(x^k-1)^{24/\sigma_1(N)}$} as automorphism of $\Lambda$. The
eta product 
\[  \eta_g(\tau) = \prod_{k|N} \eta(k \tau)^{24/\sigma_1(N)}  \]
is a cusp form for $\Gamma_0(N)$ with multiplicative coefficients. The
fixpoint lattice 
$\Lambda^g$ is the unique lattice in its genus
without roots. We lift $f_g = 1/\eta_g$ to a vector valued modular form 
\[ F_g = \sum_{M \in \Gamma_0(N)\backslash \Gamma } \, 
             f_g|_M  \, \rho_D(M^{-1}) e^0  \]
for the Weil representation $\rho_D$ of the lattice $\Lambda^g \oplus
I\!I_{1,1} \oplus \sqrt{N} I\!I_{1,1}$. Then we apply the singular theta
correspondence to $F_g$ to obtain an automorphic product $\Psi_g$
(cf.\ \cite{B3}, Theorem 13.3). We summarize this in the following
diagram 
\[ g  \mapsto  1/\eta_g  \mapsto  F_g  \mapsto  \Psi_g  \, .  \]
The automorphic form $\Psi_g$ has singular weight and is completely
reflective, i.e.\ has zeros of order one corresponding to all roots of 
$\Lambda^g \oplus I\!I_{1,1} \oplus \sqrt{N} I\!I_{1,1}$ 
(cf.\ \cite{S3}). The expansion of $\Psi_g$ at any cusp is given by 
\[
e^{\rho} \prod_{d|N} \prod_{\alpha \in (L \cap d{L'})^+} 
                        (1-e^{\alpha})^{ [1/\eta_g](-\alpha^2/2d) }   
= \sum_{ w\in W}\det(w)w(\eta_g(e^{\rho})) 
\]
where $L=\Lambda^g \oplus I\!I_{1,1}$, $\rho$ is a primitive norm $0$
vector in $I\!I_{1,1}$ and $W$ is the full reflection group of
$L$. This is the denominator identity of a generalized Kac-Moody algebra
whose real simple roots are the simple roots of $W$, i.e. the roots
$\al$ of $L$ with $(\rho,\al) = -\al^2/2$, and imaginary simple roots
are the positive multiples $n\rho$ of the Weyl vector with
multiplicity $24 \, \sigma_0((N,n))/\sigma_1(N)$. 

In the prime order case the above Lie algebras have also been
constructed by a different method in \cite{N} (cf.\ also section 14 in
\cite{B3}). 

One of the main results of \cite{S3} is that the 10 generalized Kac-Moody
algebras corresponding to the elements of squarefree order in $M_{23}$ are
the only generalized Kac-Moody algebras whose denominator identities
are completely re\-flective automorphic products of singular weight on
lattices of squarefree level splitting two hyperbolic planes. 

For $N=1$ and $2$ these Lie algebras represent the physical states of
bosonic strings moving on suitable spacetimes (cf.\ \cite{B2},
\cite{HS}). 

Let $V$ be a meromorphic vertex operator algebra of central charge
$24$ and nonzero $L_0$-eigenspace $V_1$. Schellekens \cite{ANS} proves
that either $V_1$ has dimension $24$ and $V$ is the vertex
operator algebra of the Leech lattice or $V_1$ has dimension greater
than 24. He shows that in this case there are exactly $69$ modular invariant
partition functions and describes explicitly the corresponding vector
spaces as sums of highest weight modules over affine Kac-Moody
algebras. If the monster vertex operator algebra is unique and for
each of the 69 partition functions there exists a unique vertex
operator algebra, Schellekens' result implies that there are $71$
meromorphic vertex operator algebras of central charge $24$. Up to now
these conjectures are open. 

In this paper we prove the following results. 

Let $p=2,3,5$ or $7$ and let $V$ be the prospective vertex operator
algebra $V$ in \cite{ANS} of central charge $24$, trivial fusion
algebra and spin-1 algebra $\hat{A}_{q,p}^r$ where $q=p-1$ and
$r=48/q(p+1)$. Then the character of $V$ as $\hat{A}_q^r$-module can
be written as
\[ \chi = \sum_{\gamma \in N'/N} F_{\gamma} \, \vartheta_{\gamma} \, . \]
Here $N$ is the unique lattice of minimal norm $4$ for $p=2,3,5$ and
minimal norm $6$ for $p=7$ in the genus  
\[ I\!I_{2m,0}\Big( p^{\epsilon_p (m+2)} \Big)  \]
where $m=24/(p+1)$ and 
$\epsilon_p=+$ for $p=2,5,7$ and $\epsilon_p=-$ for $p=3$, 
$\vartheta_{\gamma}$ is the theta function of $\gamma+N$ and
$F_{\gamma}$ the component corresponding to $\gamma+N$ of the lift of
$1/\eta(\tau)^m \eta(p\tau)^m$ to the lattice $N$.  

Suppose the vertex operator algebra $V$ of central charge $24$,
trivial fusion algebra and spin-1 algebra $\hat{A}_{q,p}^r$ exists and
admits a real form. Then the cohomology group of ghost number one of
the BRST-operator $Q$ acting on the vertex superalgebra
\[ V \otimes V_{I\!I_{1,1}} \otimes V_{b,c}  \]
gives a natural realization of the generalized Kac-Moody algebra
corresponding to the elements of order $p$ in $M_{23}$. 

For $p=2$ these results are also obtained in \cite{HS}. The advantage of
the methods used here is that they do not use explicit formulas for
the string functions which are in general unknown so that they can
probably be applied to all the prospective vertex operator algebras in
\cite{ANS}. 

The above result puts evidence to the conjecture that the
generalized Kac-Moody algebras whose denominator identities are
completely re\-flective automorphic products of singular weight on
lattices of squarefree level describe bosonic strings moving on
suitable spacetimes.

The paper is organized as follows.

In section $2$ we recall some results on the affine Kac-Moody algebras and
their highest weight representations.  

In section $3$ we describe the Weil representation of
$\mathit{SL}_2(\Z)$ and construct vector valued modular forms.

In section $4$ we recall some properties of vertex operator
algebras and WZW models.   

In section $5$ we show that the character of the prospective vertex
operator algebra $V$ in \cite{ANS} of spin-1 algebra $\hat{A}_{q,p}^r$
can be written in the form $\chi = \sum F_{\gamma} \,
\vartheta_{\gamma}$ as described above. 

In the last section we show that the physical states of a chiral bosonic
string with vertex algebra $V \otimes V_{I\!I_{1,1}}$ give a
realization of the generalized Kac-Moody algebra corresponding to the
elements of order $p$ in $M_{23}$. 

We thank J.\ Fuchs, G.\ H\"ohn and C.\ Schweigert for helpful
discussions and comments. We also thank the referee for suggesting
many improvements to the paper.

\section{Affine Kac-Moody algebras} \label{KMA}

In this section we recall some results on the untwisted affine Kac-Moody
algebras and their highest weight representations from \cite{K}. 

Let $g$ be a finite dimensional simple complex Lie algebra with Cartan
subalgebra $h$ and root system $\Delta$. Then $g$ has an invariant
symmetric bilinear form which is nondegenerate on $g$ and on $h$ so
that there is a natural isomorphism $\nu : h \to h^*$. We will often
write $\al^2$ for $(\al,\al)$. The coroot $\al^{\vee}$ of a root $\al$
is the inverse image of $2\al/\al^2$ under $\nu$. Let $\{
\al_1,\ldots,\al_l \}$ be a set of simple roots of $g$
and $a_{ij} = \al_j(\al^{\vee}_i)$ the corresponding Cartan matrix. We
denote the set of positive roots by $\Delta_+$. The reflections in the
hyperplanes orthogonal to the simple roots generate the Weyl group
$W$. The Lie algebra $g$ has at most two different root lengths. The
highest root 
\[ \theta = \sum_{i=1}^{l} a_i \al_i \]
is a long root and 
\[ \theta^{\vee} = \sum_{i=1}^{l} a^{\vee}_i \al^{\vee}_i \, . \]
We normalize the bilinear form such that $\theta^2 = 2$.

The untwisted affine Kac-Moody algebra corresponding to $g$ is 
\[ \hat{g} = \C[t,t^{-1}]\otimes g \oplus \C K \oplus \C d \]
where $K$ is central and
\begin{align*}
[ \, t^m\otimes x , t^n\otimes y \, ] &= 
                    t^{m+n} \otimes [x,y] + m \delta_{m+n} (x,y) K   \\
[ \, d , t^n\otimes y \, ] & = n t^n\otimes y   \, .
\end{align*}
The vector space 
\[    \hat{h} = h \oplus \C K \oplus \C d     \]
is a commutative subalgebra of $\hat{g}$. We extend a linear function $\lambda$ on $h$ to $\hat{h}$ by setting $\lambda(K) = \lambda(d) = 0$. We define linear functions $\Lambda_0$ and $\delta$ on $\hat{h}$ by
\begin{align*} 
\Lambda_0 (h \oplus \C d) &= 0 \, , & \Lambda_0 (K) &= 1 \\
\delta (h \oplus \C K)    &= 0 \, , & \delta(d)     &= 1 \, .
\end{align*}
Then 
\[    \hat{h}^* = h^* \oplus \C \Lambda_0 \oplus \C \delta     \]
and we have a natural projection 
$\hat{h}^* \to h^*, \lambda \mapsto \ov{\lambda}$ 
with $\ov{\Lambda_0} = \ov{\delta} = 0$. A linear function $\lambda$ in $\hat{h}^*$ can be written 
\[ \lambda = \ov{\lambda} + \lambda(K) \Lambda_0 + \lambda(d) \delta  \]
and $\lambda(K)$ is called the level of $\lambda$. 
We also extend the bilinear form from $g$ to $\hat{g}$ by setting
\begin{gather*}
( t^m\otimes x , t^n\otimes y ) = \delta_{m+n} (x,y) \\
( t^m\otimes x , K ) = ( t^m\otimes x , d ) = 0 \\
( K , K ) = ( d , d ) = 0 \\
(K,d) = 1 .
\end{gather*}
Define 
\[  \al_0 = \delta - \theta \, . \]
Then 
\[  \al^{\vee}_0 = K - \theta^{\vee}  \]
and $\{ \al_0, \al_1,\ldots,\al_l \}$ is a set of simple roots of
$\hat{g}$ and  $\{ \al^{\vee}_0,\al^{\vee}_1,\ldots,\al^{\vee}_l \}$
is the set of coroots. The fundamental weights $\Lambda_0, \ldots, \Lambda_l$ satisfy 
\[  \Lambda_i ( \al^{\vee}_j ) = \delta_{ij} \, ,   
    \qquad  \Lambda_i ( d ) = 0  \, . 
\]
Then $\Lambda_i = \ov{\Lambda_i} + a^{\vee}_i \Lambda_0$ and
$\ov{\Lambda_1}, \ldots, \ov{\Lambda_l}$ are the fundamental weights
of $g$. Let $b_{ij}$ be the inverse of the Cartan matrix of $g$. The
scalar products of the fundamental weights are
\[   ( \Lambda_0 , \Lambda_i ) = 0 \]
for $i=0,\ldots,l$ and
\[   ( \Lambda_i , \Lambda_j ) 
= ( \ov{\Lambda_i} , \ov{\Lambda_j} ) =  \frac{a^{\vee}_i}{a_i}b_{ij}  \]
for $i,j=1,\ldots,l$.

The Weyl vector $\rho \in \hat{h}^*$ is defined by
$\rho(\al^{\vee}_i)=1$  for $i = 0,\ldots, l$ and $\rho(d)=0$ 
(\cite{K}, p.\ 82). The level of the Weyl vector is the dual Coxeter
number $h^{\vee}$.  

We can write $\Lambda \in \hat{h}^*$ as
\[  \Lambda = \sum_{i=0}^l n_i \Lambda_i + c \delta \]
with labels $n_i = \Lambda(\al^{\vee}_i)$ and $c = \Lambda(d)$. In representation theory the value of $c$ is often unimportant. For example if two weights $\Lambda$ and $\Lambda'$ have the same labels the corresponding irreducible highest weight modules $L(\Lambda)$ and $L(\Lambda')$ are isomorphic as irreducible modules of the derived algebra $\hat{g}'=[\hat{g},\hat{g}]$. 

Let 
\[ P = \{ \lambda \in \hat{h}^* \, | \,  
          \lambda(\al^{\vee}_i) \in \Z \, \text{ for } i = 0,\ldots, l  \}
\]
be the weight lattice and 
\[ P_+ = \{ \lambda \in P \, | \,  
          \lambda(\al^{\vee}_i) \geq 0 \, \text{ for } i = 0,\ldots, l  \} 
\] 
be the set of dominant integral weights. Then
\[ P = \sum_{i=0}^{l} \Z \Lambda_i + \C \delta \, . \]
We define $P^k = \{ \lambda \in P \, | \, \lambda(K) = k   \}$ 
and $P^k_+ = P^k \cap P_+$.

Let $\Lambda \in P_+^k$ with $k>0$ and $L(\Lambda)$ the corresponding irreducible highest weight module. We define
\[ m_{\Lambda} =   \frac{(\Lambda + \rho)^2}{2(k+h^{\vee})} 
                 - \frac{\rho^2}{2h^{\vee}}                     \]
and the normalized character 
\[ \chi_{\Lambda} = e^{-m_{\Lambda}\delta} \text{ch } L(\Lambda)  \, . \]
For $\lambda  \in \hat{h}^*$ we define 
\[ m_{\Lambda, \lambda} = m_{\Lambda} - \frac{\lambda^2}{2k} \]
and the string function
\[ c_{\lambda}^{\Lambda} = e^{-m_{\Lambda,\lambda} \delta} \sum_{n \in \Z} 
   \text{mult}_{L(\Lambda)} (\lambda - n \delta) e^{-n\delta} \, .
\]
The string functions are invariant under the action of the affine Weyl
group, i.e. 
\[ c_{w(\lambda)}^{\Lambda} = c_{\lambda}^{\Lambda} \]
for $w \in \hat{W}$, so that 
\[ c_{w(\lambda) + k \gamma + a \delta}^{\Lambda} =
c_{\lambda}^{\Lambda} \]
for $w \in W, \, \gamma \in M$ and $a \in \C$ 
(\cite{K}, (12.7.9)).

The multiplicities can be calculated with Freudenthal's recursion formula
\begin{multline*}
( ( \Lambda + \rho )^2 - ( \lambda + \rho )^2 ) \,
\text{mult}_{L(\Lambda)} (\lambda) = \\
2 \sum_{\al \in \Delta_+} \sum_{j\geq 1} \,
\text{mult}_{L(\Lambda)} (\al) \,
(\lambda+ j\al,\al) \,
\text{mult}_{L(\Lambda)} (\lambda+ j\al )  \, .
\end{multline*}

Let
\[ \theta_{\lambda} = e^{k \Lambda_0} \sum_{\gamma \in M + \ov{\lambda}/k } 
                      e^{- k \gamma^2 \delta/2 + k\gamma }  \]
where $M$ is the lattice generated by the long roots of $g$ (\cite{K},
(12.7.3)). 
Then 
\[ \theta_{\lambda + k \gamma + a \delta} = \theta_{\lambda} \]
for $\gamma \in M$ and $a \in \C$ (\cite{K}, (13.2.3)).

The normalized character is given by
\[ \chi_{\Lambda} = \sum_{\lambda \in P^k \bmod (kM+\C\delta)}  
                    c_{\lambda}^{\Lambda} \theta_{\lambda}  \] 
(\cite{K}, (12.7.12)).

We obtain modular forms if we replace in the above definitions
$e^{- \delta} = q$ with $q = e^{2\pi i \tau}$. The string functions
transform under the generators of $\mathit{SL}_2(\Z)$ as 
\begin{multline*} 
c_{\lambda}^{\Lambda}(-1/\tau) = \\
|M'/kM|^{-1/2} \, (-i\tau)^{-l/2}    
\sum_{ \substack{ \Lambda' \in P^k_+ \bmod  \C \delta \\ 
                  \lambda' \in P^k \bmod  (kM+ \C \delta)   } }
S_{\Lambda,\Lambda'} \, e\big( (\ov{\lambda},\ov{\lambda'})/k \big) \,
c_{\lambda'}^{\Lambda'}(\tau) 
\end{multline*}
with
\begin{multline*}
S_{\Lambda,\Lambda'} = \\
i^{|\Delta_+|} \, |M'/(k+h^\vee)M|^{-1/2} \sum_{ w \in W}
\det(w) \, e\big( -(\ov{\Lambda}+\ov{\rho},w(\ov{\Lambda'}+ \ov{\rho} ) )/
                 (k+h^{\vee}) \big)
\end{multline*}
(\cite{K}, Theorems 13.8 and 13.10) and 
\[ c_{\lambda}^{\Lambda}(\tau+1) = 
e(m_{\Lambda,\lambda}) \, c_{\lambda}^{\Lambda}(\tau) \] 
(\cite{K}, (13.10.4)).

\section{The Weil representation} \label{mf}

In this section we describe the Weil representation of
$\mathit{SL}_2(\Z)$ and construct vector valued modular forms. 
More details can be found in \cite{S3}. 

A discriminant form is a finite abelian group $D$ with a quadratic form $D \to \Q/\Z, \, \gamma \mapsto \gamma^2/2$ such that $(\bt,\gamma) = (\bt + \gamma)^2/2 - \bt^2/2 -\gamma^2/2 \mod 1$ is a nondegenerate symmetric bilinear form. The level of $D$ is the smallest positive integer $N$ such that $N\gamma^2/2 \in \Z$ for all $\gamma \in D$. The group algebra $\C [D]$ of $D$ is the algebra with basis $\{ e^{\gamma} \, | \, \gamma \in D \}$ and products $e^{\bt} e^{\gamma} = e^{\bt + \gamma}$.

Let $L$ be an even lattice with dual lattice $L'$. Then $L'/L$ is a discriminant form with the quadratic form given by $\gamma^2/2 \! \mod 1$. Conversely every discriminant form can be obtained in this way. The signature $\sign(D) \in \Z/8\Z$ of a discriminant form is defined as the signature modulo $8$ of any even lattice with that discriminant form.

Let $D$ be a discriminant form of even signature. There is a unitary action of the group $\mathit{SL}_2(\Z)$ on $\C [D]$ defined by
\begin{align*} 
\rho_D(T) e^{\gamma}  & = e(-\gamma^2/2)\, e^{\gamma} \\
\rho_D(S) e^{\gamma}  & = \frac{e(\sign(D)/8)}{\sqrt{|D|}}
                  \sum_{\bt\in D} e((\gamma,\bt))\, e^{\bt}   
\end{align*}
where 
$S = \left( \begin{smallmatrix} 0 & -1 \\ 1 & 0 \end{smallmatrix} \right)$ and 
$T = \left( \begin{smallmatrix} 1 &  1 \\ 0 & 1 \end{smallmatrix} \right)$
are the standard generators of $\mathit{SL}_2(\Z)$. This rep\-re\-sen\-tation is called Weil representation.

Let 
\[ F(\tau) = \sum_{\gamma \in D} F_{\gamma}(\tau) e^{\gamma}  \]
be a holomorphic function on the upper halfplane with values in $\C [D]$ and $k$ an integer. Then $F$ is a modular form for $\rho_D$ of weight $k$ if 
\[ F(M\tau) = (c\tau + d)^k \rho_D(M) F(\tau)  \]
for all $M = \left( \begin{smallmatrix} a & b \\ c & d
\end{smallmatrix} \right)$ in $\mathit{SL}_2(\Z)$ and $F$ is meromorphic at $\infty$. 

We can construct modular forms for the Weil representation by lifting scalar valued modular forms on $\Gamma_0(N)$. Suppose the level of $D$ divides $N$ where $N$ is a positive integer. Let $f$ be a scalar valued modular form on $\Gamma_0(N)$ of weight $k$ and character $\chi_D$. Again we allow poles at cusps. Then
\[ F(\tau) = \sum_{M \in \Gamma_0(N)\backslash \Gamma } \, 
             f|_M (\tau) \, \rho_D(M^{-1}) \, e^0  \]
is a vector valued modular form for $\rho_D$ of weight $k$ which is invariant under the automorphisms of the discriminant form. 

Now we consider the following cases. Let $p=2,3,5$ or $7$. Then there is an
automorphism of the Leech lattice of cycle shape $1^mp^m$ with
$m=24/(p+1)$. The fixpoint lattice $\Lambda_p$ is the unique lattice
in its genus without roots. Let $I\!I_{1,1}$ be the even
unimodular Lorentzian lattice of rank $2$. Then the lattice 
\[    \Lambda_p \oplus \sqrt{p} I\!I_{1,1}    \]
has level $p$ and genus 
\[ I\!I_{2m+1,1}\Big(p^{\epsilon_p(m+2)}\Big)   \]
with $\epsilon_p=+,-,+,+$ for $p=2,3,5,7$. The eta product
\[ f(\tau) = \frac{1}{\eta(\tau)^m \eta(p\tau)^m }  \]
is a modular form for $\Gamma_0(p)$ of weight $-m$ with poles at the
cusps $0$ and $\infty$ and trivial character for $p=2,3$ and $5$ and
character $\chi \left(
\left( \begin{smallmatrix} a & b \\ c & d \end{smallmatrix} \right)
\right)
= \big( \frac{d}{7} \big)$ in the case $p=7$. We
define $T$-invariant functions $g_j$ by 
\[ f(\tau/p) = g_0(\tau) + g_1(\tau) + \ldots + g_{p-1}(\tau)  \]
with $g_j|_T(\tau) = e(j/p) g_j(\tau) $, i.e.\
\[ g_j(\tau) = \frac{1}{p} \sum_{k=0}^{p-1} e(-kj/p) 
               f\big( (\tau + k)/p \big)   \, .             \]
We lift the modular form $f$ to a modular form 
$F = \sum F_{\gamma} e^{\gamma}$ 
on $\Lambda_p \oplus \sqrt{p} I\!I_{1,1}$ using the above
construction. Then $F$ has components
\[ 
F_{\gamma}(\tau) = \begin{cases}
\, f(\tau) + g_0(\tau) \, ,  & \text{if $\gamma = 0$,} \\
\, g_j(\tau) \, ,            & \text{if $\gamma \neq 0$ and 
                               $\gamma^2/2 = -j/p \,\bmod 1$.}
\end{cases} \]
The components $F_{\gamma}$ with $\gamma^2/2 = 0 \! \mod 1$ are modular forms for $\Gamma_0(p)$ of weight $-m$ and of nontrivial quadratic character in the case $p=7$.

\section{WZW theories}

In this section we recall some results on vertex operator algebras and
WZW models. References are \cite{FHL}, \cite{H}, \cite{Fu2}, \cite{FZ}
and \cite{ANS}. 

Let $V$ be a vertex operator algebra satisfying certain regularity
properties. In particular we assume that $V$ is simple and rational,
i.e.\ $V$ is irreducible as a module over itself and $V$ has only
finitely many nonisomorphic simple modules and every module decomposes
into a finite direct sum of irreducible modules. Let $\{ M_0, \dots ,
M_n \}$ with $M_0 = V$ be the set of nonisomorphic simple modules. 
The dimensions of the intertwining spaces 
\[  N_{ij}^k = \dim \binom{M_k}{M_i \, M_j}  \] 
satisfy $N_{ij}^k =  N_{ji}^k$, $N_{0i}^j =  \delta_{i}^j$ and $N_{ij}^0 = \delta_{j i^+}$ for a unique $i^+$ depending on $i$. The fusion algebra is the free vector space $\C [I]$ where $I = \{ 0, \ldots, n \}$ with products 
\[ i \times j = \sum_{k \in I} N_{ij}^k \, k \, . \]
The fusion algebra is commutative, associative and $i \times 0 = i$. 

Zhu \cite{Z} has shown that the space of genus 1 correlation functions
of $V$ is invariant under $\mathit{SL}_2(\Z)$ and has a basis
corresponding to the modules $M_i$. We denote the $S$-matrix with
respect to this basis by $S_{ij}$. A consequence of this result is
that the characters $\chi_i$ of the modules $M_i$ are invariant under
$\mathit{SL}_2(\Z)$.

The $S$-matrix is related to the structure constants of the fusion
algebra by the Verlinde formula (cf.\ \cite{Hu})
\[  N_{ij}^k = \sum_{n \in I} \frac{S_{in}S_{jn}S_{kn}^{-1}}{S_{0n}}  \, . \]

A module $M_i$ is called a simple current if for each $j$ there is a
unique $k$ such that $i \times j = k$. Then of course $i \times i^+ =
0$. This condition is also sufficient, i.e. $M_i$ is a simple
current if and only if $N_{ii^+}^j = \delta^{j 0}$. A simple current
generates an action of a group $\Z/N\Z$ on the set of modules.  

An important problem is when there is a vertex algebra structure on a
sum over the modules $M_i$ extending $V$. In some special cases,
e.g. simple current extensions of WZW models, this problem has been
solved.  

Let $g$ be a finite dimensional simple complex Lie algebra with
affinization $\hat{g}$ and $k$ a positive integer. Then the
irreducible $\hat{g}$-module $L(k \Lambda_0)$ has a canonical
structure as a simple rational vertex operator algebra which satisfies
the above mentioned regularity properties. The simple modules are the
irreducible $\hat{g}$-modules $L(\Lambda)$ with $\Lambda \in P_+^k$ and 
$\Lambda(d) = 0$. This theory is
called the WZW model corresponding to $g$ of level $k$. Since in this
case the $S$-matrix is known the structure constants of
the fusion algebra can be calculated by the Verlinde formula. 
Fuchs \cite{Fu1} has determined the simple currents of the WZW
models. The nontrivial simple currents correspond to the weights
$\Lambda$ where $\ov{\Lambda}$ is $k$ times a cominimal weight of
$g$. If $g$ is $E_8$ there is an exceptional simple current at level
$2$. 

Let $V$ be a vertex operator algebra. Then the $L_0$-eigenspace $V_1$
is a Lie algebra under $[u,v] = u_0v$ with invariant bilinear form
$(u,v) = u_1v$. The components of the fields 
$u(z)=\sum_{m\in \Z} u_m z^{-m-1}$ and $v(z)=\sum_{n\in \Z} v_n z^{-n-1}$
satisfy the commutation relations 
\[  [u_m,v_n]= ([u,v])_{m+n} + m(u,v)\delta_{m+n}1 \, . \]
Schellekens \cite{ANS} studies the character valued partition
functions of meromorphic vertex operator algebras of central charge
$24$. With the help of the level $1$ trace identities he
proves the following result. Let $V$ be a vertex operator algebra of
central charge $24$, trivial fusion algebra and nonzero $V_1$. Then
either $\dim V_1 = 24$ and $V_1$ is commutative or $\dim V_1 > 24$ and
$V_1$ is semisimple. In the first case $V$ is the vertex operator
algebra of the Leech lattice. In the second case $V$ can be written as
a sum of modules over the affinization $\hat{V}_1$ of $V_1$. Using the
trace identities of level $2$ he shows that there are at most $69$
possibilities for the Lie algebra $V_1$ if $\dim V_1 > 24$. For each
of these possibilities he finds exactly one modular invariant
partition function and describes explicitly the decomposition of $V$
as $\hat{V}_1$-module. For many of these vector spaces it is still
open whether they have a vertex algebra structure because the
extension problem for WZW models has not been solved in general so far.
Schellekens' results suggest that there are $71$
meromorphic vertex operator algebras of central charge $24$, the
monster vertex operator algebra, the vertex operator algebra of the
Leech lattice and the $69$ vertex operator algebras with Kac-Moody symmetry.

\section{Characters of some vertex operator algebras}

In this section we show that the character of the prospective conformal
field theory in \cite{ANS} of spin-1 algebra $\hat{A}_{p-1,p}^r$ can
be written in the form $\chi = \sum F_{\gamma} \, \vartheta_{\gamma}$.

Let $p=2,3,5$ or $7$ and $q=p-1$.

First we describe some properties of the affine Kac-Moody algebra
$\hat{A}_q$. The central element $K$ is given by
\[  K = \al^{\vee}_0 + \ldots + \al^{\vee}_q  \, . \]
Let 
\[ \lambda = n_0 \Lambda_0 + \ldots + n_q \Lambda_q = (n_0,\ldots,n_q)
\]
be in $P \bmod \C \delta$. Then $\lambda$ has level 
\[ \lambda(K) = n_0 + \ldots + n_q \]
and norm
\[ (\lambda,\lambda) = \lambda^2 = \sum_{i,j=1}^q b_{ij} n_i n_j  \]
where $b_{ij}$ is the inverse of the Cartan matrix of $A_q$. 

The natural projection from $P \bmod \C \delta$ to the lattice
$A_q'$ sends 
\[   \lambda = n_0 \Lambda_0 + \ldots + n_q \Lambda_q    \]
to 
\[  \ov{\lambda} = n_1 \ov{\Lambda_1} + \ldots + n_q \ov{\Lambda_q} \,
. \]
This map induces a bijection from the weights of level $p$ in
$P \bmod \C \delta$ to the lattice $A_q'$. 

The group $A_q'/A_q \cong \Z/p\Z$ is sometimes called congruence
group and its elements congruence classes. It can be represented by
the elements $\ov{\Lambda_i}$ with 
$\ov{\Lambda_i} + \ov{\Lambda_j} = \ov{\Lambda_{i+j}}$
where addition is taken modulo $p$. 

Let $\lambda = p \Lambda_0 = (p,0,\ldots,0)$. Then the irreducible
highest weight module $L(\lambda)$ is a vertex operator algebra whose
irreducible modules correspond to the weights $(n_0,\ldots,n_q)$ where
the $n_i$ are nonnegative integers and 
$n_0 + \ldots + n_q = p$. 
The simple currents are given by the weights $p\Lambda_i$. They act by
cyclicly shifting the coefficients to the right, i.e.\ 
$(p\Lambda_i).(n_0,\ldots,n_q) = (n_{q+1-i},\ldots,n_{q-i})$, and form a
group isomorphic to $\Z/p\Z$. 

We describe some properties of the simple currents.

Let $\lambda$ be a weight of level $p$ and $s$ a simple current. Then
$\ov{\lambda}$ and $\ov{s.\lambda}$ are in the same class. 

Furthermore

\begin{prp} \label{invofstringfcts}
The string functions of $\hat{A}_q$ of level $p$ are invariant
under the following action of a simple current $s$
\[ c_{\lambda}^{\Lambda} = c_{s.\lambda}^{s.\Lambda} \, . \]
\end{prp}
{\em Proof:} There is a diagram automorphism $\phi$ acting on
$\bigoplus_{i=0}^l \C \al_i$ such that 
$\phi.\al_i \bmod \C\delta = s.(\al_i \bmod \C\delta)$ for
$i=0,\dots,l$. We extend $\phi$ to $\hat{h}^*$ by $\phi.\Lambda_0 =
s.\Lambda_0$. Then there is a unique map $\phi$ on $\hat{h}$ satisfying 
$(\phi.\lambda)(\phi.y) = \lambda(y)$ for all
$\lambda \in \hat{h}^*$ and $y \in  \hat{h}$. 
This map gives an isomorphic realization
and therefore an isomorphism of $\hat{g}$. Hence
\[ c_{s.\lambda}^{s.\Lambda} = 
c_{\phi.\lambda + c \delta}^{\phi.\Lambda + C \delta} =
c_{\phi.\lambda}^{\phi.\Lambda} =
c_{\lambda}^{\Lambda} \, . \]
\eop\\

The isomorphism from the discriminant group to the group of simple
currents sending $\ov{\Lambda_i}$ to $p\Lambda_i$ satisfies

\begin{prp} \label{corrclasscurr}
Let $\mu$ be a weight of level $0$ such that 
$\ov{\mu} = \ov{\Lambda_i} \! \mod A_q$. Then 
\[ c^\Lambda_{\lambda+p\mu} = c^\Lambda_{(p{\Lambda}_i). \lambda} \]
for all $\Lambda$ in $P^p_+$ and all $\lambda$ in $P^p$.
\end{prp}
{\em Proof:} Let $\lambda = (n_0,\ldots,n_q)$ be a weight of level $p$
and $\mu = (m_0,\ldots,m_q)$ \vspace{0.1em} a weight of level $0$ such
that 
$\ov{\mu} = \ov{\Lambda_1} \! \mod A_q$. Then
\begin{align*}
w_1 w_2 \ldots w_q (\lambda + p \mu)
&= (n_q, n_0, \ldots , n_{q-1}) \mod A_q \\
&= (p\Lambda_1). \lambda \mod A_q 
\end{align*}
so that
\[ c_{\lambda + p \mu}^{\Lambda} =
   c_{w_1 w_2 \ldots w_q (\lambda + p \mu)}^{\Lambda} =
   c_{(p\Lambda_1). \lambda}^{\Lambda} \, . 
\]
This implies the statement. \eop\\

Let $V$ be the prospective vertex operator algebra in \cite{ANS} of 
spin-1 algebra $\hat{A}_{q,p}^r$ where $r=48/q(p+1)$. Then $V$ is a
sum of irreducible highest weight modules of $\hat{A}_q^r$, the weight
of each factor $\hat{A}_q$ having level $p$. We denote the set of
highest weight vectors by $M$. 
The set $M$ is invariant under the natural action of a 
subgroup $G$ of the group of simple currents and decomposes into
$G$-orbits $G\backslash M$. In the appendix we list $G$ and orbit
representatives of $G\backslash M$ together with their multiplicities.

Using the isomorphism from the discriminant group to the group of
simple currents we can consider $G$ as subgroup of
${A_q^r}'/A_q^r$. We denote by $(A_q^r,G)$ the rational lattice obtained by
gluing the elements of $G$ to $A_q^r$ and analogously by
$(A_q^r,G^{\perp})$ the rational lattice obtained by gluing the
elements of the orthogonal complement $G^{\perp}$ to
$A_q^r$. Note that $(A_q^r,G^{\perp})$ is the dual lattice of
$(A_q^r,G)$. 

The character of $V$ as $\hat{A}_q^r$-module is given by

\begin{align*}
\chi_{V} 
&= \sum_{\Lambda \in M} \mult(\Lambda) \, \, \chi_{\Lambda} \\
&= \sum_{\Lambda \in M} \mult(\Lambda) \, \prod_{i=1}^r \, 
   \sum_{\lambda_i \in A_q'/pA_q} \,
   c_{\lambda_i}^{\Lambda_i}(\tau) \, \theta_{\lambda_i}(\tau,z_i,u_i)
   \\
&= \sum_{\lambda \in {A_q^r}' /pA_q^r} \,  
   \sum_{\Lambda \in M} \mult(\Lambda) \, 
   \prod_{i=1}^r \,
   c_{\lambda_i}^{\Lambda_i}(\tau) \, \theta_{\lambda_i}(\tau,z_i,u_i)
   \, .
\end{align*}
We have

\begin{prp} \label{nonzerocontrib}
The only nonzero contributions in the above expression for the character of $V$ come from weights $\lambda$ in $(A_q^r,G^{\perp})$. 
\end{prp}
{\em Proof:} A case-by-case analysis shows $M\subset
(A_q^r,G^{\perp})$. It is easy to see that $c_{\lambda}^{\Lambda}=0$
if $\Lambda$ and $\lambda$ are not in the same class. Therefore any
nonzero contribution to the character of $V$ comes from a weight
$\lambda$ in $(A_q^r,G^{\perp})$. \eop\\

It follows

\begin{prp} \label{formulaf}
The character of $V$ as $\hat{A}_q^r$-module can be written as 
\[  \chi_{V} = \sum_{\lambda \in N'/N} 
\tilde{F}_{\lambda} \, \vartheta_{\lambda}         \]
where $N=\sqrt{p}(A_q^r,G)$, $\vartheta_{\lambda}$ is the theta function of the coset $\lambda+N$ and 
\[ \tilde{F}_{\lambda} =
   \sum_{\Lambda \in G\backslash M} \, \sum_{g \in G} \, 
   \frac{\mult( \Lambda)}{|G_{\Lambda}|}  
   \prod_{i=1}^r c_{g_i.\sqrt{p}\lambda_i}^{\Lambda_i}    \, .            \]
\end{prp}
{\em Proof:} By Prop\-o\-si\-tions \ref{nonzerocontrib} and \ref{invofstringfcts} we have
\begin{align*}
\chi_{V} 
&= \sum_{\lambda \in {A_q^r}' /pA_q^r} \,  
   \sum_{\Lambda \in M} \mult(\Lambda) \, 
   \prod_{i=1}^r \,
   c_{\lambda_i}^{\Lambda_i}(\tau) \, \theta_{\lambda_i}(\tau,z_i,u_i)
   \\
&= \sum_{\lambda \in (A_q^r,G^{\perp})/pA_q^r } \,  
   \sum_{\Lambda \in M} \mult(\Lambda) \, 
   \prod_{i=1}^r \,
   c_{\lambda_i}^{\Lambda_i}(\tau) \, \theta_{\lambda_i}(\tau,z_i,u_i)
   \\
&= \sum_{\lambda \in (A_q^r,G^{\perp})/pA_q^r } \, 
   \sum_{\Gamma \in G\backslash M} \, \sum_{\Lambda \in G.\Gamma}
   \mult(\Gamma) \, \prod_{i=1}^r \,
   c_{\lambda_i}^{\Lambda_i}(\tau) \, \theta_{\lambda_i}(\tau,z_i,u_i) \\
&= \sum_{\lambda \in (A_q^r,G^{\perp})/pA_q^r } \,
   \sum_{\Lambda \in G\backslash M } \, \sum_{ g \in G } 
   \frac{\mult( \Lambda)}{|G_{\Lambda}|} \, \prod_{i=1}^r \,
   c_{\lambda_i}^{g_i.\Lambda_i}(\tau) \, \theta_{\lambda_i}(\tau,z_i,u_i) \\
&= \sum_{\lambda \in (A_q^r,G^{\perp})/pA_q^r } \,
   \sum_{\Lambda \in G\backslash M } \, \sum_{ g \in G } 
   \frac{\mult( \Lambda)}{|G_{\Lambda}|}
   \, \prod_{i=1}^r \,
   c_{g_i.\lambda_i}^{\Lambda_i}(\tau) \, \theta_{\lambda_i}(\tau,z_i,u_i) \\
&= \sum_{\lambda \in (A_q^r,G^{\perp})/pA_q^r } \,
   \prod_{i=1}^r \, \theta_{\lambda_i}(\tau,z_i,u_i) \,
   \sum_{\Lambda \in G\backslash M } \, \sum_{ g \in G } 
   \frac{\mult( \Lambda)}{|G_{\Lambda}|}
   \, \prod_{i=1}^r \, c_{g_i.\lambda_i}^{\Lambda_i}(\tau)  
\end{align*}
where $G_{\Lambda}$ denotes the stabilizer of $\Lambda$ in $G$. Using Proposition \ref{corrclasscurr} we get 
\begin{align*}
\chi_{V}
&= \sum_{\lambda \in (A_q^r,G^{\perp})/p(A_q^r,G) } \,
   \sum_{ \mu \in p(A_q^r,G)/pA_q^r } \,
   \prod_{i=1}^r \, \theta_{\lambda_i+\mu_i}(\tau,z_i,u_i)  \\
&  \qquad \qquad \qquad \cdot
   \sum_{\Lambda \in G\backslash M } \, \sum_{ g \in G } 
   \frac{\mult( \Lambda)}{|G_{\Lambda}|}
   \, \prod_{i=1}^r \, c_{g_i.(\lambda_i + \mu_i)}^{\Lambda_i}(\tau)   \\
&= \sum_{\lambda \in (A_q^r,G^{\perp})/p(A_q^r,G) } \,
   \sum_{ \mu \in p(A_q^r,G)/pA_q^r } \,
   \prod_{i=1}^r \, \theta_{\lambda_i+\mu_i}(\tau,z_i,u_i)  \\
&  \qquad \qquad \qquad \cdot
   \sum_{\Lambda \in G\backslash M } \, \sum_{ g \in G } 
   \frac{\mult( \Lambda)}{|G_{\Lambda}|}
   \, \prod_{i=1}^r \, c_{s_i.g_i.\lambda_i}^{\Lambda_i}(\tau)   \\
&= \sum_{\lambda \in (A_q^r,G^{\perp})/p(A_q^r,G) } \,
   \sum_{ \mu \in p(A_q^r,G)/pA_q^r } \,
   \prod_{i=1}^r \, \theta_{\lambda_i+\mu_i}(\tau,z_i,u_i)  \\
&  \qquad \qquad \qquad \cdot
   \sum_{\Lambda \in G\backslash M } \, \sum_{ g \in G } 
   \frac{\mult( \Lambda)}{|G_{\Lambda}|}
   \, \prod_{i=1}^r \, c_{g_i.\lambda_i}^{\Lambda_i}(\tau) \, .  
\end{align*} 
Here $s_i$ is the simple current of $\hat{A}_q$ such that 
$c_{g_i.\lambda_i + g_i.\mu_i}^{\Lambda_i} 
= c_{s_i.g_i.\lambda_i}^{\Lambda_i}$. 
\vspace{0.2em} 
Note that $s=(s_1, \ldots, s_r)$
is in $G$ because $\mu=(\mu_1, \ldots, \mu_r)$ is in
$p(A_q^r,G)/pA_q^r$. Now 
\begin{eqnarray*}
\lefteqn{\sum_{ \mu \in p(A_q^r,G)/pA_q^r } \,
\prod_{i=1}^r \, \theta_{\lambda_i+\mu_i}(\tau,z_i,u_i)   }          \\
&=& 
\sum_{ \mu \in p(A_q^r,G)/pA_q^r } \,
\sum_{ \nu \in pA_q^r } \, 
\prod_{i=1}^r \, e(pu_i) \, e\big( \tau (\lambda_i+\mu_i+\nu_i)^2/2p 
                             + (\lambda_i+\mu_i+\nu_i,z_i) \big)  \\ 
&=& 
\sum_{ \mu \in p(A_q^r,G)} \,
\prod_{i=1}^r \, e(pu_i) \, e\big( \tau (\lambda_i+\mu_i)^2/2p 
                             + (\lambda_i +\mu_i,z_i) \big)        \\
&=& 
\sum_{ \mu \in \sqrt{p}(A_q^r,G)} \,
\prod_{i=1}^r \, e(pu_i) \, e\big( \tau (\lambda_i+\sqrt{p}\mu_i)^2/2p 
                             + (\lambda_i +\sqrt{p}\mu_i,z_i) \big)  \\
&=& 
\sum_{ \mu \in \sqrt{p}(A_q^r,G)} \,
\prod_{i=1}^r \, e(pu_i) \, e\big( \tau (\lambda_i/\sqrt{p} +\mu_i)^2/2 
                             + (\lambda_i/\sqrt{p} +\mu_i,\sqrt{p}z_i) \big) 
\end{eqnarray*}
so that 
\begin{align*}
\chi_{V} 
&= \sum_{\lambda \in (A_q^r,G^{\perp})/p(A_q^r,G) } \,
   \sum_{ \mu \in \sqrt{p}(A_q^r,G)} \\
&  \qquad \qquad 
   \prod_{i=1}^r \, e(pu_i) \, e\big( \tau (\lambda_i/\sqrt{p} +\mu_i)^2/2 
                             + (\lambda_i/\sqrt{p} +\mu_i,\sqrt{p}z_i) \big) \\
&  \qquad \qquad \cdot 
   \sum_{\Lambda \in G\backslash M } \, \sum_{ g \in G } 
   \frac{\mult( \Lambda)}{|G_{\Lambda}|}
   \, \prod_{i=1}^r \, c_{g_i.\lambda_i}^{\Lambda_i}(\tau)                 \\
&= \sum_{\lambda \in (1/\sqrt{p})(A_q^r,G^{\perp})/ 
                                    \sqrt{p}(A_q^r,G) } \,
   \sum_{ \mu \in \sqrt{p}(A_q^r,G)} \\
&  \qquad \qquad
   \prod_{i=1}^r \, e(pu_i) \, e\big( \tau (\lambda_i +\mu_i)^2/2 
                             + (\lambda_i +\mu_i,\sqrt{p}z_i) \big) \\
&  \qquad \qquad \cdot 
   \sum_{\Lambda \in G\backslash M } \, \sum_{ g \in G } 
   \frac{\mult( \Lambda)}{|G_{\Lambda}|}
   \, \prod_{i=1}^r \, c_{g_i.\sqrt{p}\lambda_i}^{\Lambda_i}(\tau) \, . 
\end{align*}
This implies the proposition. \eop\\

Note that the functions $\tilde{F}_{\gamma}$ are invariant under $G$, 
the permutations $\Sym(G)$ of the $r$ components which leave $G$ invariant and by Proposition \ref{invofstringfcts} also under
$\hat{W}^r$.  

\begin{prp}
The lattice $N=\sqrt{p}(A_q^r,G)$ has genus 
\[ I\!I_{2m,0}\Big( p^{\epsilon_p (m+2)} \Big)  \]
where $m=24/(p+1)$ and $\epsilon_p=+,-,+,+$ for $p=2,3,5,7$. In the
case $p=2$ the 2-adic Jordan components are even. The minimal norm of
$N$ is 4 for $p=2$, $3$ and $5$, and $6$ for $p=7$. This is the largest
possible minimal norm of a lattice in this genus and $N$ is the unique lattice up to isomorphism with this minimal norm.  
\end{prp}
{\em Proof:} It is easy to see that $N$ is an even lattice. Furthermore we have $pN' = \sqrt{p}(A_q^r,G^{\perp}) \subset \sqrt{p}(A_q^r,G) = N$ because $G^{\perp}\subset G$. Using $|G|^2=p^{24/q-2}$ we get
\begin{align*}
\det (N) 
&= \det (A_q^r,G)(p) \\
&= p^{rq} \det (A_q^r,G)\\
&= p^{rq} \det (A_q^r) / |G|^2 \\
&= p^{rp}/ |G|^2 \\
&= p^{m+2} \, . 
\end{align*}
This implies that $N$ has the stated genus (cf.\ \cite{CS}, p.\ 386
ff., Theorem 13). The minimal norm of $N$ follows from the minimal
distance of $G$ considered as linear code in $\F_p^r$. We find that
the minimal norm of $N$ is $4$ for $p=2,3$ and $5$, and
$6$ for $p=7$. We leave the proof of the other statements to the
reader. \eop\\

Now we determine the function 
\[  \tilde{F} = \sum_{\gamma \in D} \tilde{F}_{\gamma} e^{\gamma}\, ,
\] 
where $D$ is the discriminant form of $N$, explicitly. 

\begin{thm} 
The function $\tilde{F}$ is a modular form of weight $-24/(p+1)$ for the Weil representation of $N$.
\end{thm}
{\em Proof:} We have to show that
\[ \tilde{F}_{\gamma} (\tau+1) = e(-\gamma^2/2) \, \tilde{F}_{\gamma} (\tau) \]
and 
\[ \tilde{F}_{\gamma} (-1/\tau) = 
\frac{e(\sign(D)/8)}{\sqrt{|D|}} \, \tau^k
\sum_{\bt\in D} e\big( (\gamma,\bt) \big) \, \tilde{F}_{\bt} (\tau)   \]
with $k=-24/(p+1)$ for all $\gamma$ in $D$. 

To prove these equations we proceed as follows. 
We choose a set of functions $\{ \tilde{F}_{\gamma_1},
\tilde{F}_{\gamma_2}, \ldots, \tilde{F}_{\gamma_n} \}$ such that each
$\tilde{F}_{\gamma}$ is conjugate to exactly one
$\tilde{F}_{\gamma_j}$ under the action of $G$, $\Sym(G)$ and
$\hat{W}^r$.
 
Using the $T$-invariance of the string functions we verify that  
\[ \tilde{F}_{\gamma_j} (T\tau) 
   = e(-\gamma_j^2/2) \, \tilde{F}_{\gamma_j} (\tau) \, . \]
This implies that the $\tilde{F}_{\gamma}$ transform correctly under
$T$.  

Let $\gamma \in D$. We define constants $c_{\gamma, \gamma_l}$ by 
\[ \sum_{\bt\in D} e\big( (\gamma,\bt) \big) \, \tilde{F}_{\bt} 
 = \sum_{l} c_{\gamma,\gamma_l} \tilde{F}_{\gamma_l}  \, . \]
If $\gamma$ is equivalent to $\gamma_j$ under the above symmetries
then 
\[ c_{\gamma,\gamma_l} = c_{\gamma_j,\gamma_l} \, . \]
The action of $S$ on $\tilde{F}_{\gamma_j}$ is given by 
\[ \tilde{F}_{\gamma_j}(S\tau) = 
\sum_{\Lambda \in G\backslash M} \, \sum_{g \in G}  
\frac{\mult( \Lambda)}{|G_{\Lambda}|}  
\prod_{i=1}^r c_{g_i.\sqrt{p}(\gamma_j)_i}^{\Lambda_i}(S\tau) \, . \]
We determine the $S$-matrix of the string functions by computer
calculations using the formula in section 2 and write
$\tilde{F}_{\gamma_j}(S\tau)$ as polynomial in the string
functions. We check that 
\[ \tilde{F}_{\gamma_j}(S\tau) = \frac{e(\sign(D)/8)}{\sqrt{|D|}} \,
\tau^k \sum_{l} c_{\gamma_j,\gamma_l} \tilde{F}_{\gamma_l} (\tau) \, . \] 
To see this it is helpful to replace the weight $\lambda$ in the
string function $c_{\lambda}^{\Lambda}$ in the expression of
$\tilde{F}_{\gamma_j}(S\tau)$ and of the $\tilde{F}_{\gamma_l}$ by the
unique dominant weight in the $\hat{W}$-orbit of $\lambda$ in the set
of weights of $L(\Lambda)$. 
This shows that $\tilde{F}_{\gamma_j}$ and the $\tilde{F}_{\gamma}$
have the desired transformation behaviour under $S$. \eop \\

Since the theta function of a lattice transforms under the dual Weil representation of the corresponding lattice it is now obvious that the character 
\[  \chi_{V} = \sum_{\gamma \in N'/N} 
\tilde{F}_{\gamma} \, \vartheta_{\gamma}         \]
is invariant under $\mathit{SL}_2(\Z)$. 

As in section \ref{mf} let 
\[   f(\tau) = \frac{1}{\eta(\tau)^m \eta(p\tau)^m }   \]
and define $T$-invariant functions $g_j$ by 
\[ f(\tau/p) = g_0(\tau) + g_1(\tau) + \ldots + g_{p-1}(\tau)  \]
where $g_j|_T(\tau) = e(j/p) g_j(\tau)$. Then

\begin{thm} \label{compf}
The modular form $\tilde{F}$ is equal to the lift $F$ of the modular
form $f$ to $N$.
\end{thm}
{\em Proof:} Since the lattices $N$ and 
$\Lambda_p \oplus \sqrt{p} I\!I_{1,1}$ have the same signature
modulo $8$ and isomorphic discriminant forms the components of $F =
\sum F_{\gamma} e^{\gamma}$ are given by 
\[ 
F_{\gamma}  = \begin{cases}
\, f  + g_0 \, , & \text{if $\gamma = 0$,} \\
\, g_j      \, , & \text{if $\gamma \neq 0$ and 
                               $\gamma^2/2 = -j/p \,\bmod 1$.}
\end{cases} \]
We calculate the first nonvanishing coefficient of the string functions
by computer using Freudenthal's formula and determine the singular
coefficients of $\tilde{F}$. It turns out that they
are equal to the singular coefficients of $F$. Hence the difference of
$\tilde{F}$ and $F$ is a holomorphic modular form of
negative weight which is finite at $\infty$ and therefore must be
$0$. \eop\\

If explicit formulas for the string functions are available then
Theorem \ref{compf} can also be proved in the following way. The
modular properties of the string functions imply that $\tilde{F}_{\gamma}$ in
the formula of Proposition \ref{formulaf} is a modular form of weight
$k$ and some level $N$ with poles at cusps. The same is true for
$F_{\gamma}$. Hence we can deduce equality of these functions by
comparing sufficiently many coefficients. For $\hat{A}_1$ and
$\hat{A}_2$ there are explicit formulas for the string functions
determined by Kac and Peterson \cite{KP}. We have used them to verify
the statement of Theorem \ref{compf} in the cases $p=2$ and $3$.

Unfortunately explicit formulas for the string functions are known
only in a very few cases. The advantage of the proof of Theorem
\ref{compf} given above is that it only needs the first nonvanishing
coefficient of the string functions which is easy to determine using Freudenthal's formula. Therefore our method can also be
applied to the other prospective vertex operator algebras in
\cite{ANS}.

\section{Construction of some generalized Kac-Moody algebras as bosonic strings}
We show that the physical states of a chiral bosonic
string with vertex algebra $V \otimes V_{I\!I_{1,1}}$ give a
realization of the generalized Kac-Moody algebra corresponding to the
elements of order $p$ in $M_{23}$.

We assume now that the prospective vertex operator algebra $V$ in
\cite{ANS} with spin-1 algebra $\hat{A}_{p-1,p}^r$ exists and has a
real form. 

We will work in this section over the real numbers.

Let $V_{I\!I_{1,1}}$ be the vertex algebra of the Lorentzian lattice
$I\!I_{1,1}$. Recall that the $b,\!c\,$-ghost system of the bosonic
string is described by the vertex superalgebra of the lattice $\Z$. It
carries a conformal structure of weight $-26$. We call the
corresponding vertex operator superalgebra $V_{b,c}$. There is an
action of the BRST-operator $Q$ with $Q^2=0$ on the vertex superalgebra
\[ V \otimes V_{I\!I_{1,1}} \otimes V_{b,c}\, . \]
The cohomology group of ghost number one has a Lie bracket \cite{LZ}
and we denote this Lie algebra by $G$. The vertex algebra $V$ is
graded by the rational lattice $N'$ so that $G$ is graded by 
$N' \oplus I\!I_{1,1}$. The no-ghost theorem (cf.\ Theorem 5.1 in
\cite{B2}) implies that the graded dimensions are given by $\dim
G_{\al} = 2m+2$ if $\al =0$ and $\dim G_{\al} = [F_{\al}](-\al^2/2)$
if $\al \neq 0$. Here $F_{\al}$ is the component of $F$ corresponding
to $\al \! \mod I\!I_{1,1}$ (cf.\ Theorem \ref{compf}). Furthermore
Theorem 2 in \cite{B4} shows that $G$ is a generalized Kac-Moody
algebra. 

In order to obtain an even grading lattice $L$ we rescale $N'\oplus
I\!I_{1,1}$ by $p$. It is easy to see that $L$ has genus
$I\!I_{2m+1,1}( p^{+m} )$. The lattice $\Lambda_p \oplus
I\!I_{1,1}$ has the same genus. It follows from Eichler's theory of
spinor genera that there is only one class in this genus so that $L$
is isomorphic to $\Lambda_p \oplus I\!I_{1,1}$. Thus we have

\begin{prp}
The Lie algebra $G$ is a generalized Kac-Moody algebra graded by the
Lorentzian lattice $L$. The Cartan subalgebra has dimension $2m+2$ and 
\[ 
\dim \, G_{\al}  = \begin{cases}
\,\, [f](-\al^2/2) \, , & \text{if $\al \in L\backslash pL'$}, \\
\,\, [f](-\al^2/2) + [f](-\al^2/2p) \, , & \text{if $\al \in pL'$},
\end{cases} \]
for nonzero $\al$, where  
\[ f(\tau) = \frac{1}{\eta(\tau)^m \eta(p\tau)^m } 
           = q^{-1} + m + \ldots \, .   \]
\end{prp}

We recall that $[f](n)$ denotes the coefficient at $q^n$ in the Fourier
expansion of the function $f$. 

Let $\rho$ be a primitive norm $0$ vector in $I\!I_{1,1}$. Then $\rho$
is a Weyl vector for the reflection group $W$ of $L = \Lambda_p \oplus
I\!I_{1,1}$ (cf.\ \cite{B1}). 

\begin{thm}
The denominator identity of $G$ is 
\begin{multline*}
e^{\rho} 
\prod_{\al \in L^+} (1-e^{\al})^{ [f](-\al^2/2) }
 \prod_{\al \in p{L'}^+} (1-e^{\al})^{ [f](-\al^2/2p) } \\
= \sum_{ w\in W}\det(w) \, w \!
\left( e^{\rho} \prod_{n=1}^{\infty} (1-e^{n\rho})^m (1-e^{pn\rho})^m 
\right) \, .  
\end{multline*}
The real simple roots of $G$ are the simple roots of $W$, i.e. the norm $2$ vectors in $L$ with $(\rho,\al) = -1$ and the norm $2p$ vectors in $pL'$ with $(\rho,\al) = -p$. The imaginary simple roots of $G$ are the positive multiples $n\rho$ of the Weyl vector with multiplicity $m$ if $p \! \not| \,n$ and multiplicity $2m$ if $p|n$.
\end{thm}
{\em Proof:} We only have to prove the second statement. Let $K$ be
the generalized Kac-Moody algebra with root lattice $L$, Cartan
subalgebra $L\otimes \R$ and simple roots as stated in the theorem. We
lift $f$ to a vector valued modular form $F$ on $L \oplus \sqrt{p}
I\!I_{1,1}$. Note that $F$ admits the same description as the lift of
$f$ on $\Lambda_p \oplus \sqrt{p} I\!I_{1,1}$ in section
\ref{mf}. Then we apply the singular theta correspondence to $F$ to
obtain an automorphic form $\Psi$ of singular weight. The expansion of
$\Psi$ at any cusp is given by
\begin{multline*}
e^{\rho} 
\prod_{\al \in L^+} (1-e^{\al})^{ [f](-\al^2/2) }
 \prod_{\al \in p{L'}^+} (1-e^{\al})^{ [f](-\al^2/2p) } \\
= \sum_{ w\in W}\det(w) \, w \!
\left( e^{\rho} \prod_{n=1}^{\infty} (1-e^{n\rho})^m (1-e^{pn\rho})^m 
\right) \, .  
\end{multline*}
This is the denominator identity of $K$. We see that $G$ and $K$ have
the same root multiplicities. The product in the denominator identity
determines the simple roots of $G$ because we have fixed a Cartan
subalgebra and a fundamental Weyl chamber. It follows that $G$ and $K$
have the same simple roots and are isomorphic. \eop

\section*{Appendix}

Below we list the groups $G$ as linear codes in $\F_p^r$ and
orbit representatives of $G\backslash M$ together with their
multiplicities. 

In the case $p=2$ the glue code $G$ is the binary Hamming code of
length $16$ and orbit representatives are 
\begin{itemize}
\item $(2,0)^{16}$
\item $8\times (1,1)^{16}$
\item The remaining orbit representatives can be described as
  follows. In the dual binary Hamming code of length $16$, for every
  codeword of weight $8$, identify the $1$-components with the highest
  weight $(1,1)$ and for the $0$-components allow all combinations of
  $(2,0)$ and $(0,2)$ such that both of these highest weights appear
  an odd number of times.
\end{itemize}

\vspace*{0.9em}
In the case $p=3$ the glue code $G$ is the ternary zero sum 
code of length $6$ and orbit representatives are
\begin{itemize}
\item $(3,0,0)^6$
\item $(1,1,1)^4(3,0,0)^2$ and all permutations
\item $(2,0,1)^5(0,1,2)$ and $(2,1,0)^5(0,2,1)$
\item $6\times(1,1,1)^6$.
\end{itemize} 

\vspace*{0.9em}
In the case $p=5$ the glue code is $\F_5^{ 2}$ and orbit
representatives are
\begin{itemize}
\item $(5,0,0,0,0)^2$
\item $(2,0,1,0,2)^2$
\item $(2,0,0,2,1)(3,0,1,1,0)$ and $(3,0,1,1,0)(2,0,0,2,1)$
\item $(1,1,1,1,1)(1,0,0,1,3)$ and $(1,0,0,1,3)(1,1,1,1,1)$
\item $4\times(1,1,1,1,1)^2$.
\end{itemize}

\vspace*{0.9em}
In the case $p=7$ the glue code is $\F_7$ and orbit representatives are
\begin{itemize}
\item $(7,0,0,0,0,0,0)$
\item $(2,0,0,1,3,0,1)$ and $(2,1,0,3,1,0,0)$
\item $(2,0,0,2,0,3,0)$
\item $(1,0,1,0,1,2,2)$
\item $3\times(1,1,1,1,1,1,1)$.
\end{itemize}

\vspace*{10mm}
\noindent
Thomas Creutzig,\\
DESY Theory Group, Notkestrasse 85, Bldg.\ 2a, 22603 Hamburg,
Germany, thomas.creutzig@desy.de\\[3mm]
Alexander Klauer,\\
University of Mannheim, Department of Mathematics, 68131 Mannheim,
Germany, 
aklauer@rumms.uni-mannheim.de \\[3mm]
Nils R.\ Scheithauer,\\
University of Edinburgh, Maxwell Institute for Mathematical Sciences,
School of Mathematics, James Clerk Maxwell Building,  
Mayfield Road, Edin\-burgh \mbox{EH9 3JZ}, United Kingdom, Nils.Scheithauer@ed.ac.uk

\end{document}